\documentclass[11pt]{article}

\usepackage{graphicx}

\def\be{\begin{equation}}
\def\ee{\end{equation}}
\newcommand{\ff}[1]{{\mbox{\boldmath $#1$}}}
\def\a{\alpha}
\def\b{\beta}
\def\x{\ff{x}}

\title{Engineering Optimisation by Cuckoo Search }

\author{Xin-She Yang \\
Department of Engineering \\
University of Cambridge \\
Trumpington Street \\
Cambridge CB2 1PZ, UK \\
\and
Suash Deb \\
Dept of Computer Science \& Engineering \\
C. V. Raman College of Engineering \\
Bidyanagar, Mahura, Janla \\
Bhubaneswar 752054, INDIA }

\date{}
\begin{document}
\maketitle

\begin{abstract}

A new metaheuristic optimisation algorithm, called Cuckoo Search (CS), was developed recently by
Yang and Deb (2009). This paper presents a more extensive comparison study using some standard
test functions and newly designed stochastic test functions. We then apply the CS algorithm
to solve engineering design optimisation problems, including the design of springs
and welded beam structures. The optimal solutions obtained by CS are far better
than the best solutions obtained by an efficient particle swarm optimiser.
We will discuss the unique search features used in CS and
the implications for further research. \\ \\

\noindent {\bf Reference to this paper} should be made as follows: \\
Yang, X.-S., and Deb, S. (2010), ``Engineering Optimisation by Cuckoo Search'',
{\it Int. J. Mathematical Modelling and Numerical Optimisation},
Vol.~{\bf 1}, No.~4, 330--343 (2010).

\end{abstract}


\section{Introduction}

Most design optimisation problems in engineering are often highly nonlinear, involving
many different design variables under complex constraints. These constraints can be written
either as simple bounds such as the ranges of material properties, or as nonlinear relationships
including maximum stress, maximum deflection, minimum load capacity, and geometrical configuration.
Such nonlinearity often results in multimodal response landscape. Subsequently, local search
algorithms such as hill-climbing and Nelder-Mead downhill simplex methods are not suitable,
only global algorithms should be used so as to obtain optimal solutions (Deb 1995, Arora 1989,
Yang 2005, Yang 2008).

Modern metaheuristic algorithms have been developed with an aim to carry out
global search, typical examples are genetic algorithms (Glodberg 1989),
particle swarm optimisation (PSO) (Kennedy and Eberhart 1995, Kennedy et al 2001).
The efficiency of metaheuristic algorithms can be attributed to the
fact that they imitate the best features in nature, especially the selection of the fittest
in biological systems which have evolved by natural selection over millions of years.
Two important characteristics of metaheuristics are:
intensification and diversification (Blum and Roli 2003, Gazi and Passino 2004, Yang 2009).
Intensification intends to search around the current best solutions and select the best
candidates or solutions, while diversification makes
sure that the algorithm can explore the search space more efficiently, often by randomization.

Recently, a new metaheuristic search algorithm, called Cuckoo Search (CS),
has been developed by Yang and Deb (2009). Preliminary studies show that
it is very promising and could outperform existing algorithms such as PSO.
In this paper, we will further study CS and validate it against test functions
including stochastic test functions. Then, we will apply it to solve design optimisation problems
in engineering. Finally, we will discuss the unique features of
Cuckoo Search and propose topics for further studies.

\section{Cuckoo Search}

In order to describe the Cuckoo Search more clearly, let us briefly review the
interesting breed behaviour of certain cuckoo species. Then, we will outline
the basic ideas and steps of the proposed algorithm.

\subsection{Cuckoo Breeding Behaviour}

Cuckoo are fascinating birds, not only because of the beautiful sounds they can make,
but also because of their aggressive reproduction strategy. Some species such as the
{\it ani} and {\it Guira} cuckoos lay their eggs in communal nests, though
they may remove others'
eggs to increase the hatching probability of their own eggs (Payne et al 2005).
Quite a number of species engage the obligate brood parasitism by laying their eggs in the
nests of other host birds (often other species). There are three basic types of
brood parasitism: intraspecific brood parasitism, cooperative breeding, and
nest takeover. Some host birds can engage direct conflict with the intruding
cuckoos. If a host bird discovers
the eggs are not its owns, it will either throw these alien eggs away
or simply abandons its nest and builds  a new nest elsewhere. Some cuckoo species
such as the New World brood-parasitic {\it Tapera}
have evolved in such a way that female parasitic cuckoos are often
very specialized in the mimicry in colour and pattern of the eggs of a few chosen
host species (Payne et al 2005). This reduces the probability of their eggs
being abandoned and thus increases their reproductivity.

Furthermore, the timing of egg-laying of some species
is also amazing. Parasitic cuckoos often choose a nest where the
host bird just laid its own eggs. In general,
the cuckoo eggs hatch slightly earlier than their host eggs. Once
the first cuckoo chick is hatched, the first instinct action it
will take is to evict the host eggs by blindly propelling
the eggs out of the nest, which increases the cuckoo chick's
share of food provided by its host bird (Payne et al 2005).
Studies also show that a cuckoo chick can
also mimic the call of host chicks to gain access to more feeding opportunity.

\subsection{L\'evy Flights}

In nature, animals search for food in a random or quasi-random manner.
In general, the foraging path of an animal is effectively a random walk because the next move is
based on the current location/state and the transition probability to
the next location. Which direction it chooses depends implicitly on a
probability which can be modelled mathematically. For example, various studies
have shown that the flight behaviour of many animals and insects has demonstrated
the typical characteristics of L\'evy flights (Brown et al 2007, Reynods and Frye 2007,
Pavlyukevich 2007).

A recent study by Reynolds and Frye (2007) shows that
fruit flies or {\it Drosophila melanogaster}, explore their landscape using a series
of straight flight paths punctuated by a sudden $90^{o}$ turn, leading to
a L\'evy-flight-style intermittent scale-free search pattern.
Studies on human behaviour such as the Ju/'hoansi hunter-gatherer foraging
patterns also show the typical feature of L\'evy flights.
Even light can be related to L\'evy flights (Barthelemy et al 2008).
Subsequently, such behaviour has been applied to
optimization and optimal search, and preliminary results show its
promising capability (Shlesinger 2006, Pavlyukevich 2007).

\subsection{Cuckoo Search}

For simplicity in
describing our new Cuckoo Search (Yang and Deb 2009), we now use the following
three idealized rules:
\begin{itemize}
\item Each cuckoo lays one egg at a time, and dumps it in a randomly chosen nest;

\item  The best nests with high quality of eggs (solutions) will carry over to the next generations;

\item The number of available host nests is fixed, and a host can discover an alien egg
    with a probability $p_a \in [0,1]$. In this case,
   the host bird can either throw the egg  away or abandon the nest so as to
   build a completely new nest in a new location.
\end{itemize}

For simplicity, this last assumption can be approximated by a fraction $p_a$ of the $n$ nests
being replaced by new nests (with new random solutions at new locations).
For a  maximization problem, the quality or fitness of a solution can simply be proportional
to the objective function. Other forms of fitness can be defined in a similar
way to the fitness function in genetic algorithms.

Based on these three rules, the basic steps of the Cuckoo Search (CS)
can be summarised as the pseudo code shown in Fig. \ref{kuk-fig-100}.

\begin{figure}
\begin{center}
\begin{minipage}[c]{0.9\textwidth}
\hrule \vspace{5pt}
\indent \quad Objective function $f(\x), \;\; \x=(x_1, ..., x_d)^T$; \\
\indent \quad Initial a population of $n$ host nests $\x_i \; (i=1,2,...,n)$; \\
\indent \quad {\bf while} ($t<$MaxGeneration) or (stop criterion); \\
\indent \qquad Get a cuckoo (say $i$) randomly by L\'evy flights; \\
\indent \qquad Evaluate its quality/fitness $F_i$;  \\
\indent \qquad Choose a nest among $n$ (say $j$) randomly; \\
\indent \qquad  {\bf if } ($F_i>F_j$), \\
\indent \qquad \qquad Replace $j$ by the new solution;\\
\indent \qquad  {\bf end} \\
\indent \qquad  Abandon a fraction ($p_a$) of worse nests \\
\indent \qquad \qquad [and build new ones at new locations via L\'evy flights]; \\
\indent \qquad Keep the best solutions (or nests with quality solutions); \\
\indent \qquad Rank the solutions and find the current best;  \\
\indent \quad {\bf end while} \\
\indent \quad Postprocess results and visualisation;
\hrule \vspace{5pt}
\caption{Cuckoo Search (CS). \label{kuk-fig-100} }
\end{minipage}
\end{center}
\end{figure}

When generating new solutions $\x^{(t+1)}$ for, say cuckoo $i$,
a L\'evy flight is performed
\be \x^{(t+1)}_i=\x_i^{(t)} + \a \oplus \textrm{L\'evy}(\lambda), \ee
where $\a>0$ is the step size which should be related to the scales of the problem of
interest. In most cases,  we can use $\a=O(1)$.
The product $\oplus$ means entry-wise multiplications.
L\'evy flights  essentially provide a random walk while their random steps
are drawn from a L\'evy distribution for large steps
\be \textrm{L\'evy} \sim u = t^{-\lambda}, \quad (1 < \lambda \le 3), \ee
which has an infinite variance with an infinite mean.
Here the consecutive jumps/steps of a cuckoo  essentially form a random walk
process which obeys a power-law step-length distribution with a heavy tail.

It is worth pointing out that, in the real world, if a cuckoo's egg is 
very similar to a host's eggs, then this cuckoo's egg is less likely to be discovered, 
thus the fitness should be related to the difference in solutions.  
Therefore, it is a good idea to do a random walk in a biased way with 
some random step sizes. A demo version is attached in the Appendix 
(this demo is not published in the actual paper, but as a supplement to help
readers to implement the cuckoo search correctly).

\section{Implementation and Validation}

\subsection{Validation and Parameter Studies}

\begin{figure}
 \centerline{\includegraphics[height=3.5in,width=4in]{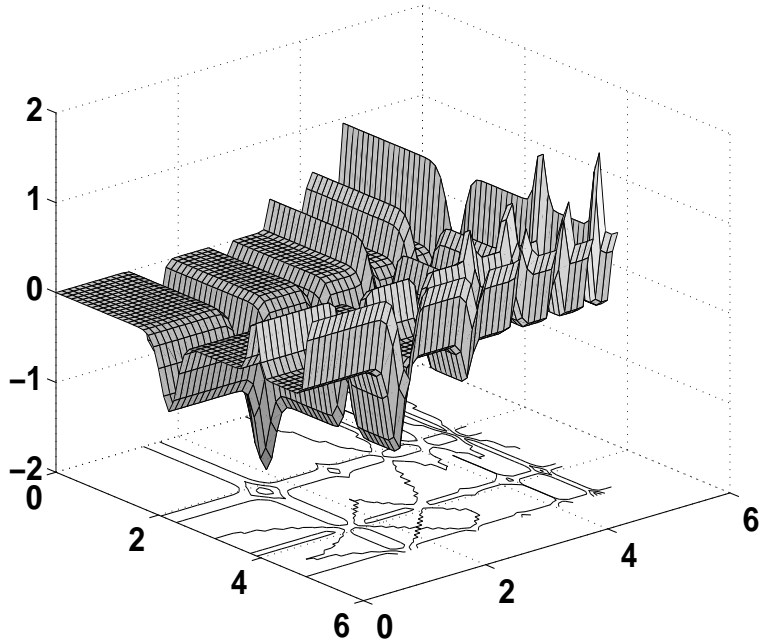} }
\caption{The landscape of Michaelwicz's 2D function. \label{kuk-fig-150} }
\end{figure}

It is relatively easy to implement the algorithm, and then we have to
benchmark it using test functions with analytical
or known solutions. There are many benchmark test functions and there is no standard list or collection,
though extensive descriptions of various functions do exist in literature
(Floudas et al 1999,  Hedar 2005, Molga and Smutnicki 2005). For example,
Michalewicz's test function has many local optima
\be f(\x)=-\sum_{i=1}^d \sin (x_i) \Big[\sin (\frac{i x_i^2}{\pi}) \Big]^{2m}, \;\; (m=10), \ee
in the domain $0 \le x_i \le \pi$ for $i=1,2,...,d$ where $d$ is the number of dimensions.
The global mimimum $f_*\approx -1.801$ occurs at $(2.20319,1.57049)$ for $d=2$,
while $f_* \approx -4.6877$ for $d=5$. In the 2D case, its 3D landscape
is shown Fig. \ref{kuk-fig-150}.

\begin{figure}
 \centerline{\includegraphics[height=2.25in,width=2.5in]{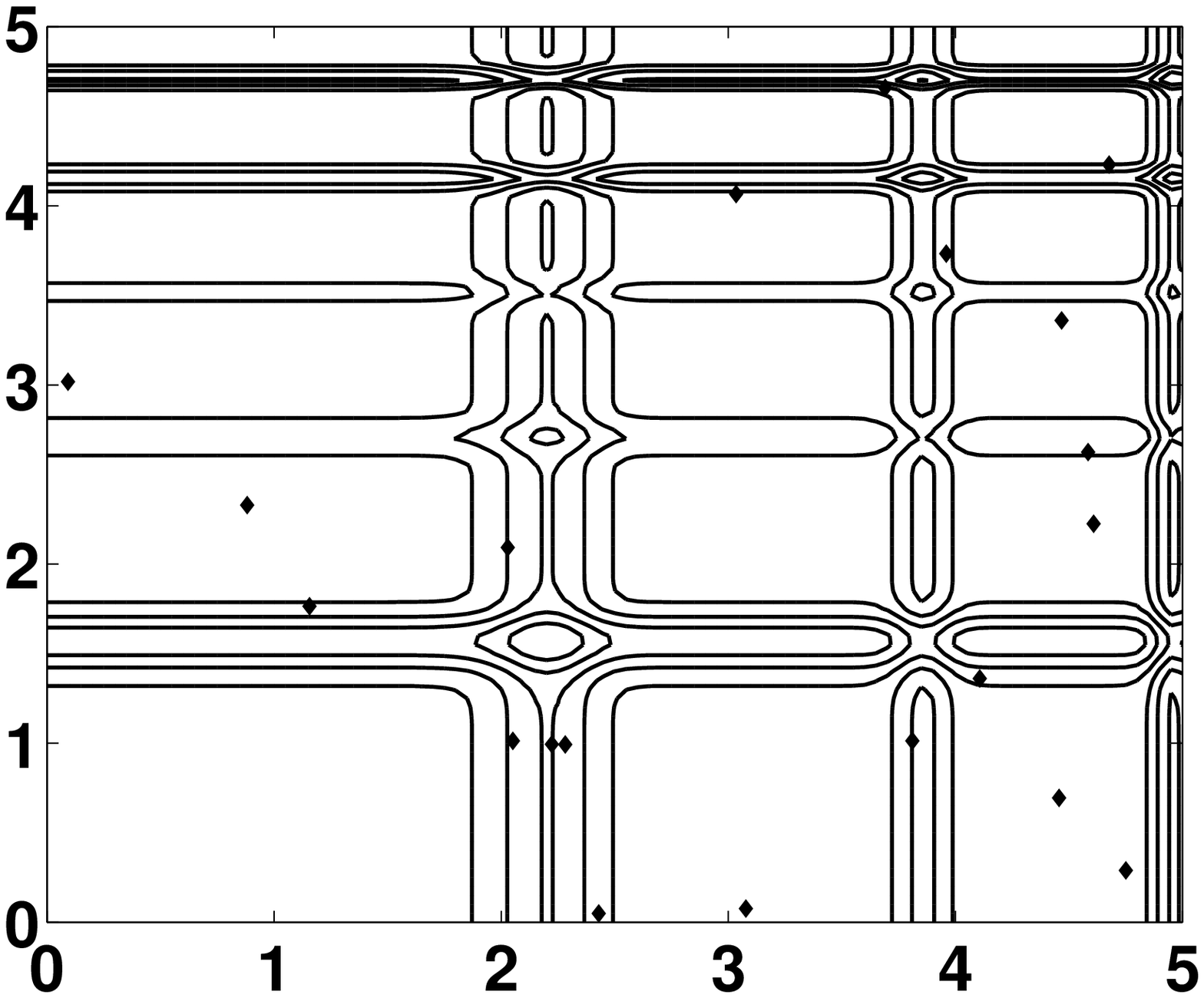}
 \includegraphics[height=2.25in,width=2.5in]{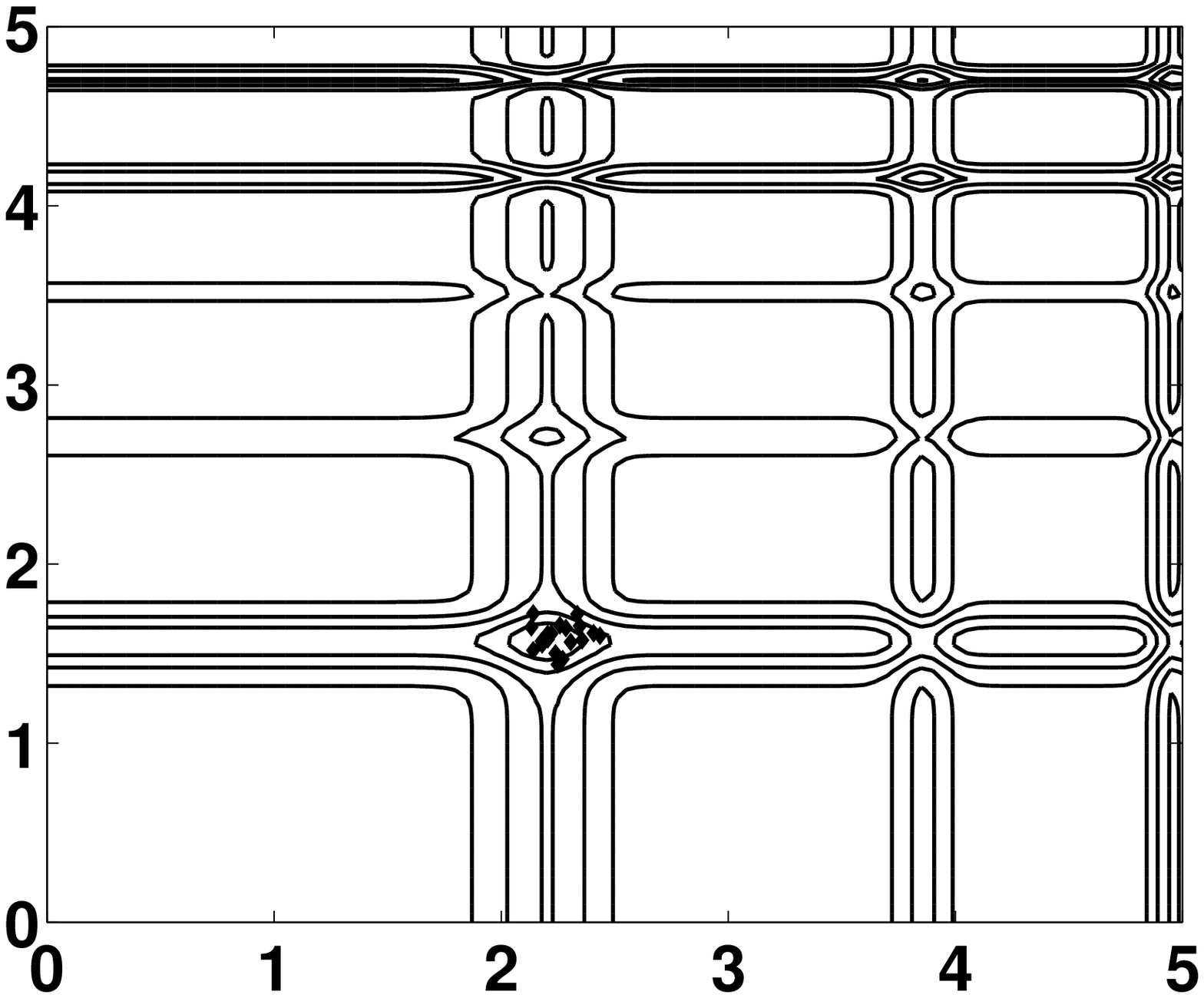}
}
\caption{Initial locations of 20 nests in Cuckoo Search, and their final
locations are marked with $\diamond$. \label{kuk-fig-200} }
\end{figure}

The global optimum in 2D can easily be
found using Cuckoo Search, and the results are shown in Fig. \ref{kuk-fig-200}
where the final locations of the nests are marked with $\diamond$. Here we have
used $n=20$ nests, $\a=1$ and $p_a=0.25$.  From the figure, we can see
that, as the optimum is approaching, most nests aggregate towards the global optimum.
In various simulations, we also notice that nests are also distributed at
different (local) optima in the case of multimodal functions.
This means that CS can find all the optima simultaneously
if the number of nests are much higher than the number of local optima.
This advantage may become more significant when dealing with multimodal
and multiobjective optimization problems.

We have also tried to vary the number of host nests (or the population size $n$)
and the probability $p_a$. We have used $n=5,10,15,20,50, 100,150, 250, 500$
and $p_a=0, 0.01, 0.05, 0.1, 0.15, 0.2, 0.25$, $0.4, 0.5$.
From our simulations, we found that $n=15$ to $25$ and $p_a=0.15$ to $0.30$ are sufficient for
most optimization problems. Results and analysis also imply that the convergence rate, to
some extent, is not sensitive to the parameters used. This means that the fine adjustment
of algorithm-dependent parameters  is not needed for any given problems.
Therefore, we will use $n=20$ and $p_a=0.25$ in the
rest of the simulations, especially for the comparison studies presented later.

\subsection{Standard Test Functions}

Various test functions in literature are designed to
test the performance of optimization algorithms (Chattopadhyay 1971, Schoen 1993,
Shang and Qiu 2006). Any new optimization algorithm should also be
validated and tested against these benchmark functions.  In our simulations, we have used the
following test functions.

De Jong's first function is essentially a sphere function
\be f(\x)=\sum_{i=1}^d x_i^2, \quad x_i \in [-5.12,5.12], \ee
whose global minimum $f(\x_*)=0$ occurs at $\x_*=(0,0,...,0)$.
Here $d$ is the dimension.

The generalized Rosenbrock's function is given by
\be f(\x)=\sum_{i=1}^{d-1} \Big[ (1-x_i)^2 + 100 (x_{i+1}-x_i^2)^2 \Big], \ee
which has a unique global minimum $f_*=0$ at $\x_*=(1,1,...,1)$.

Schwefel's test function is multimodal
\be f(\x)=\sum_{i=1}^d \Big[ -x_i \sin (\sqrt{|x_i|} \;) \Big],
\;\; -500 \le x_i \le 500, \ee
whose global minimum  $f_*=-418.9829d$ is at $x_i^*=420.9687 (i=1,2,...,d)$.

Ackley's function is also multimodal
\be f(\x)=-20 \exp\Bigg[-0.2 \sqrt{\frac{1}{d} \sum_{i=1}^d x_i^2} \; \Bigg]
-\exp[\frac{1}{d} \sum_{i=1}^d \cos (2 \pi x_i)] + (20+e), \ee
with the global minimum $f_*=0$ at $\x_*=(0,0,...,0)$ in the
range of $-32.768 \le x_i \le 32.768$ where $i=1,2,...,d$.

Rastrigin's test function
\be f(\x)=10 d + \sum_{i=1}^d [x_i^2 - 10 \cos (2 \pi x_i) ], \ee
has a unique global minimum $f_*=0$ at $(0,0,...,0)$ in a hypercube
$-5.12 \le x_i \le 5.12$ where $i=1,2,...,d$.

Easom's test function has a sharp tip
\be f(x,y)=-\cos(x) \cos(y) \exp[-(x-\pi)^2 - (y-\pi)^2], \ee
in the domain $(x,y) \in [-100,100] \times [-100,100]$.
It has a global minimum of $f_*=-1$ at $(\pi,\pi)$ in a very small region.

Griewangk's test function has many local minima
\be f(\x)=\frac{1}{4000} \sum_{i=1}^d x_i^2 - \prod_{i=1}^d \cos(\frac{x_i}{\sqrt{i}})+1, \ee
but a unique global mimimum $f_*=0$ at $(0,0,...,0)$ for all $-600 \le x_i \le 600$
where $i=1,2,...,d$.

\subsection{Stochastic Test Functions}

Almost all the test functions in literature are deterministic. It is usually
more difficult for algorithms to deal with stochastic functions. We have designed
some stochastic test functions for such a purpose.

\begin{figure}
\centerline{\includegraphics[height=2.25in,width=2.5in]{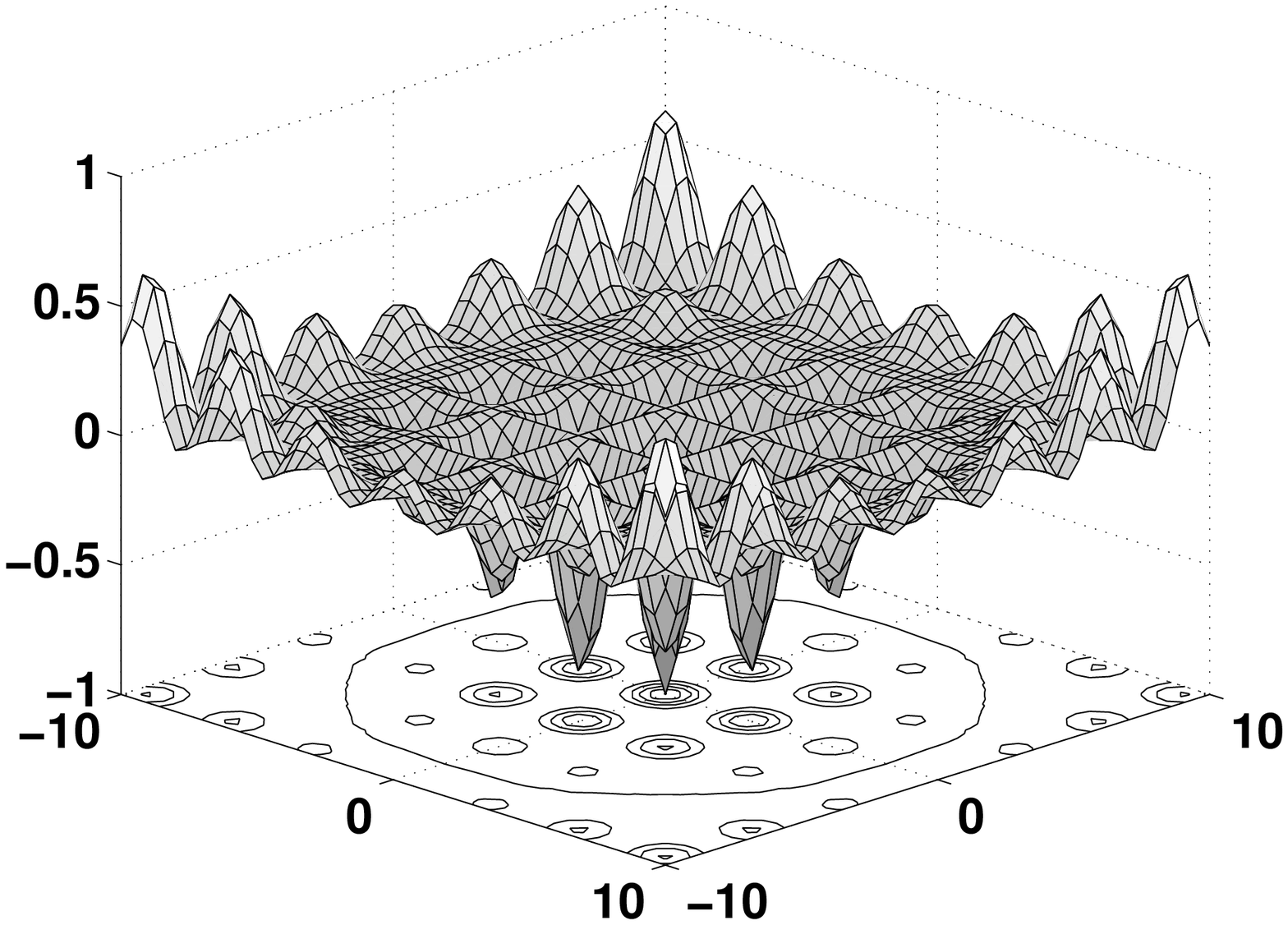}
\includegraphics[height=2.25in,width=2.5in]{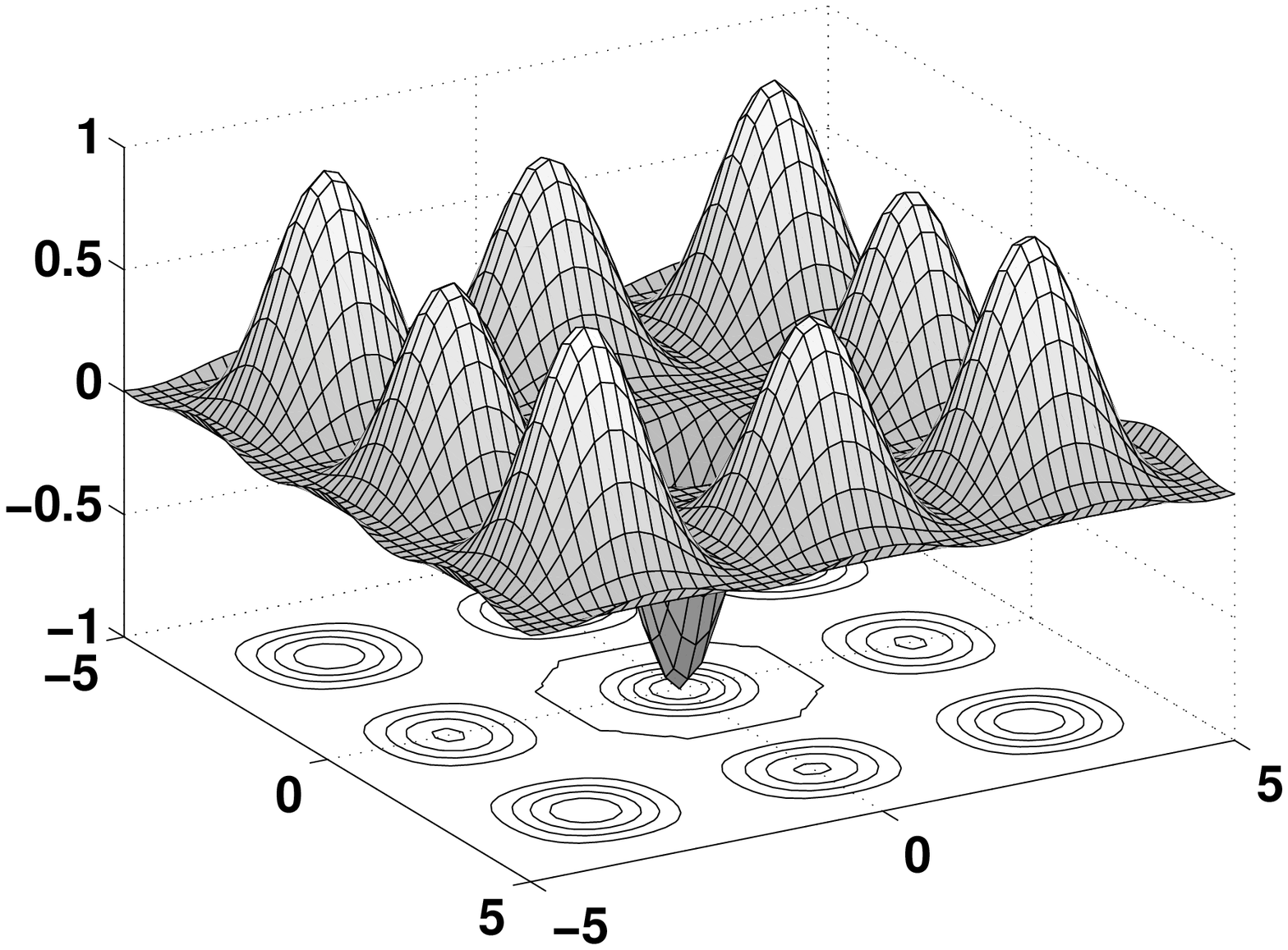} }
\caption{Landscape of stochastic function (\ref{equ-yang-50}) for small $\epsilon$
(left) and large $\epsilon$ (right). \label{fig-3ab} }
\end{figure}

The first test function designed by Yang (2010)
looks like a standing-wave function with a region of defects
\be f(\x)=\Big[ e^{-\sum_{i=1}^d (x_i/\b)^{2m}} - 2 e^{-\sum_{i=1}^d \epsilon_i (x_i-\pi)^2} \Big]
\cdot \prod_{i=1}^d \cos^2 x_i, \quad m=5, \label{equ-yang-50} \ee
which has many local minima and the unique global minimum $f_*=-1$ at $\x_*=(\pi,\pi,...,\pi)$
for $\b=15$ within the domain $-20 \le x_i \le 20$ for $i=1,2,...,d$.
Here the random variables $\epsilon_i \; (i=1,2,...,d)$ are uniformly distributed in $(0,1)$.
For example, if all $\epsilon_i$ are relatively small (say order of $0.05$), a snapshot
of the landscape in 2D is shown in Fig. \ref{fig-3ab}, while for higher values such as $0.5$
the landscape is different, also shown in Fig. \ref{fig-3ab}.

Yang's second test function is also multimodal but it has a singularity
\be f(\x) =\Big( \sum_{i=1}^d \epsilon_i |x_i| \Big) \exp \Big[ -\sum_{i=1}^d \sin (x_i^2) \Big],
\label{fun-equ-100} \ee
which has a unique global minimum $f_*=0$ at $\x_*=(0,0,...,0)$
in the domain $-2 \pi \le x_i \le 2 \pi$ where $i=1,2,...,d$ (Yang 2010).
This function is singular at the optimum $(0,...,0)$. Similarly,
$\epsilon_i$ should be drawn from a uniform distribution in $[0,1]$ or Unif[0,1].
In fact, using the same methodology, we can turn many determistic functions into stochastic test
functions. For example, we can extend Robsenbrock's function as
the following stochastic function
\be f(\x)=\sum_{i=1}^{d-1} \Big[(1-x_i)^2 + 100 \epsilon_i (x_{i+1}-x_i^2)^2  \Big], \ee
where $\epsilon_i$ should be drawn from Unif[0,1].
Similarly, we can also extend De Jong's function into its corresponding stochastic form
\be f(\x) =\sum_{i=1}^d \epsilon_i x_i^2, \ee
which still has the same global minimum $f_*=0$ at $(0,0,...,0)$, despite its stochastic nature
due to the factor $\epsilon_i$.
For stochastic functions, most deterministic algorithms such as hill climbing
and Nelder-Mead downhill simplex method would simply fail. However, we can see later
that most metaheuristic algorithms such as PSO and CS are still robust.

\subsection{Simulations and Comparison}

Recent studies indicate that PSO can
outperform genetic algorithms (GA) and
other conventional algorithms (Goldberg 1989, Kennedy et al 2001, Yang 2008).
This can be attributed partly to the broadcasting ability of the current
best estimates, potentially leading to a better and quicker convergence rate towards the
optimality.  A general framework for evaluating statistical
performance of evolutionary algorithms
has been discussed in detail by Shilane et al (2008).

Now we can compare the Cuckoo Search with PSO and genetic
algorithms for various  test functions.
After implementing these algorithms using
Matlab, we have carried out extensive simulations and each algorithm has been
run at least 100 times so as to carry out meaningful statistical analysis.
The algorithms stop when the variations of function values are less than
a given tolerance $\epsilon \le 10^{-5}$.
The results are summarised in  Table 1 where the numbers are in the
format: average number of evaluations $\pm$ one standard deviation (success rate), so
$3321 \pm 519 (100\%)$ means that the average number (mean) of function
evaluations is 3321 with a standard deviation of 519. The success rate
of finding the global optima for this algorithm is $100\%$.
The functions used in the Table are (1) Michaelwicz ($d=16$),
(2) Rosenrbrock ($d=16$), (3) De Jong ($d=32$), (4) Schwefel ($d=32$),
(5) Ackley ($d=128$), (6) Rastrigin, (7) Easom, (8) Griewangk,
(9) Yang's first stochastic function, (10) Yang's second stochastic function,
(11) Generalised Robsenbrock's function with stochastic components, and
(12) De Jong's stochastic function.

\begin{figure}
\centerline{\includegraphics[height=2in,width=2.5in]{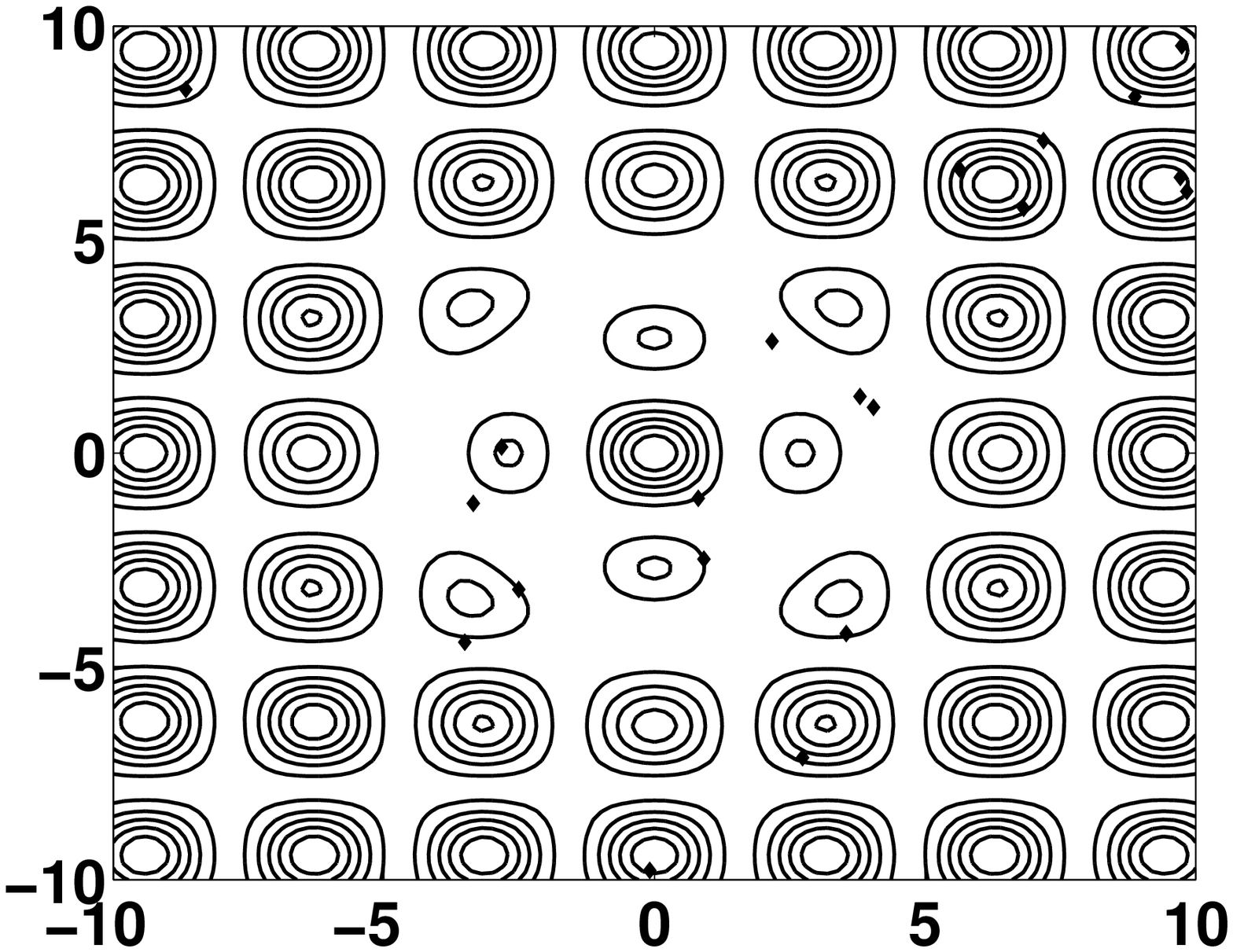}
\includegraphics[height=2in,width=2.5in]{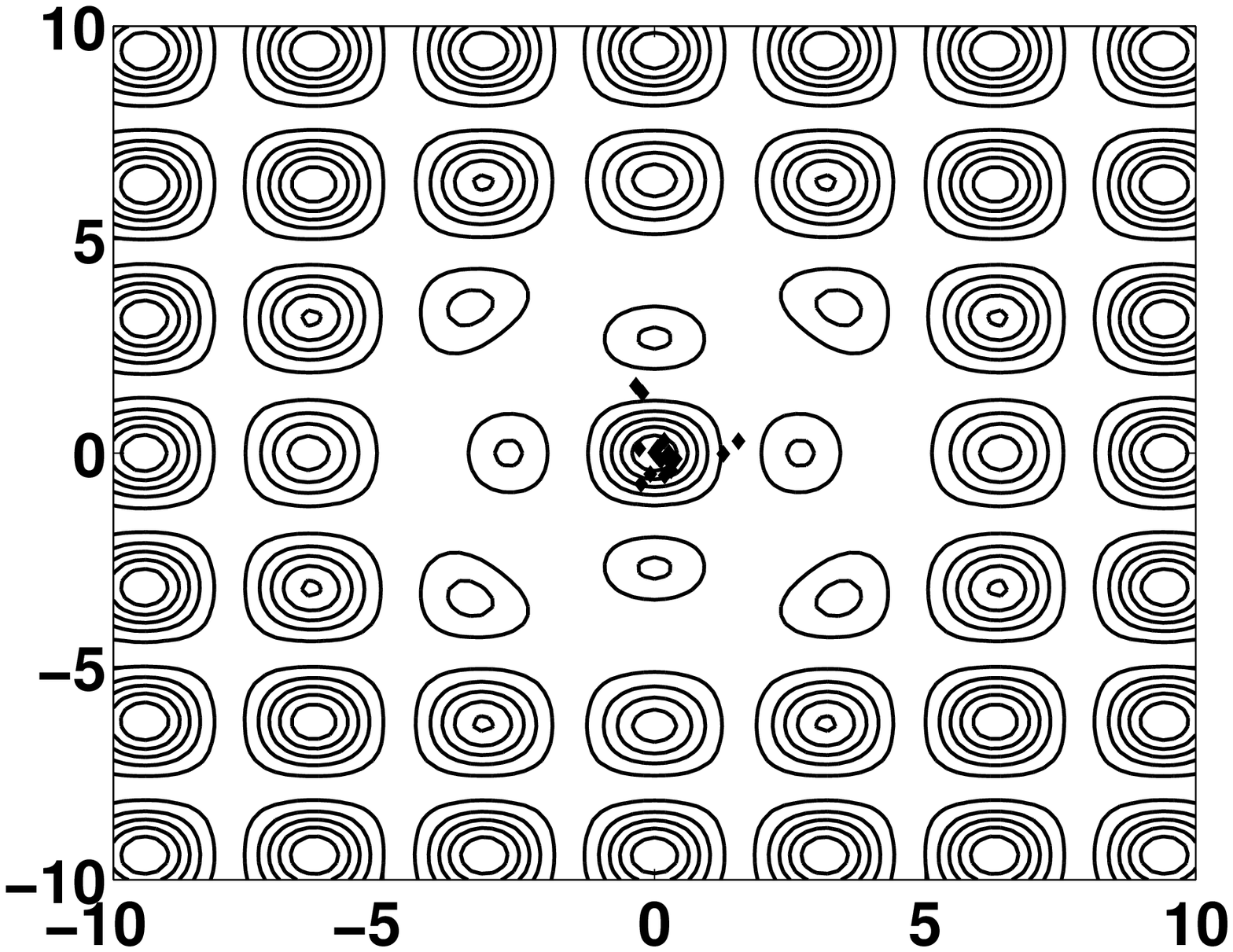} }
\caption{The initial location of 20 nests (left) for function
(\ref{equ-yang-50}) and their final locations after $15$ iterations (right).
\label{fig-4ab} }
\end{figure}

\begin{table}[ht]
\caption{Comparison of CS with genetic algorithms and particle swarm optimisation}
\centering
\begin{tabular}{ccccc}
\hline \hline
Functions & GA & PSO & CS \\
\hline
 (1)   & $89325 \pm 7914 (95 \%)$   & $6922 \pm 537 (98\%)$ & $3221 \pm 519 (100\%)$ \\
 (2) & $55723 \pm 8901 (90\%)$ & $32756 \pm 5325 (98\%)$ & $5923 \pm 1937 (100\%) $ \\
 (3)  & $15232 \pm 1270 (100\%)$  & $10079 \pm 970 (100\%)$ & $3015 \pm 540 (100\%)$\\
 (4) & $23790 \pm 6523 (95\%)$ & $92411 \pm 1163 (97\%)$  & $4710 \pm 592 (100\%)$ \\
 (5) & $32720 \pm 3327 (90\%)$ & $23407 \pm 4325 (92\%)$ & $4936 \pm 903 (100\%)$ \\
 (6) & $110523 \pm 5199 (77 \%)$ & $79491 \pm 3715 (90\%)$ & $10354 \pm 3755 (100\%)$ \\
 (7) & $19239 \pm 3307 (92\%)$  & $17273 \pm 2929 (90\%)$ & $6751 \pm 1902 (100\%)$ \\
 (8) & $70925 \pm 7652 (90\%)$  & $55970 \pm 4223 (92\%)$ & $10912 \pm 4050 (100\%)$ \\
 (9) & $79025 \pm 6312 (49\%)$  & $34056 \pm 4470 (90\%)$ & $11254 \pm 2733 (99\%)$ \\
 (10)& $35072 \pm 3730 (54\%)$  & $22360 \pm 2649 (92\%)$ & $8669 \pm 3480 (98\%)$ \\
 (11)& $63268 \pm 5091 (40\%)$  & $49152 \pm 6505 (89\%)$ & $10564 \pm 4297 (99\%)$ \\
 (12)& $24164 \pm 4923 (68\%)$  & $11780 \pm 4912 (94\%)$ & $7723 \pm 2504 (100\%)$ \\
\hline
\end{tabular}
\end{table}

We can see that CS is much more efficient in finding the global optima
with higher success rates. Each function evaluation is virtually instantaneous
on a modern personal computer. For example, the computing time for 10,000 evaluations
on a 3GHz desktop is about 5 seconds. In addition, for stochastic functions,
genetic algorithms do not perform well, while PSO is better. However, CS is
far more promising.

\section{Engineering Design}

Design optimisation is an integrated part of designing
any new products in engineering and industry. Most design problems are
complex and multiobjective, sometimes even the optimal solutions of interest do not
exist. In order to see how the CS algorithm may perform, we now use two standard
but well-known test problems.

\subsection{Spring Design Optimisation}

Tensional and/or compressional springs are used widely in engineering.
A  standard spring design problem has three design variables: the wire diameter $w$,
the mean coil diameter $d$, and the length (or number of coils) $L$.

The objective is to minimise the weight of the spring, subject to various
constraints such as maximum shear stress, minimum deflection, and
geometrical limits. For detailed description, please refer to earlier
studies (Belegundu 1982, Arora 1989, Cagnina et al 2008).
This problem can be written compactly as
\be \textrm{Minimise } f(\x)=(L+2) w^2 d, \qquad \qquad  \ee
subject to
\be \begin{array}{lll}
g_1(\x) =1-\frac{d^3 L}{71785 w^4} \le 0, \\ \\
g_2(\x) =1-\frac{140.45 w}{d^2 L} \le 0, \\ \\
g_3(\x) =\frac{2(w + d)}{3}-1 \le 0, \\ \\
g_4(\x)= \frac{d (4 d-w)}{w^3 (12566  d - w)} + \frac{1}{5108 w^2} -1 \le 0,
\end{array} \ee
with the following limits
\be 0.05 \le w \le 2.0, \quad 0.25 \le d \le 1.3, \quad 2.0 \le L \le 15.0. \ee

Using Cuckoo Search, we have obtained the same or slightly better 
solutions than the best solution obtained by Cagnina et al (2008)
\be f_*=0.012665  \quad \textrm{ at } \;\; (0.051690, 0.356750, 11.287126), \ee
but cuckoo search uses significantly fewer evaluations. 

\subsection{Welded Beam Design}

The so-called  welded beam design is another standard test
problem for constrained design optimisation  (Ragsdell and Phillips 1976,
Cagnina et al 2008). The problem has four design variables: the width $w$
and length $L$ of the welded area, the depth $h$ and thickness $h$ of the main
beam. The objective is to minimise the overall fabrication cost, under
the appropriate constraints of shear stress $\tau$,
bending stress $\sigma$, buckling load $P$ and maximum end deflection $\delta$.

The problem can be written as
\be \textrm{minimise } \; f(\x)=1.10471 w^2 L + 0.04811 d h (14.0+L), \ee
subject to
\be
\begin{array}{lll}
g_1(\x)=w -h \le 0, \vspace{2pt} \\ \vspace{3pt}
g_2(\x) =\delta(\x) - 0.25 \le 0, \\ \vspace{3pt}
g_3(\x)=\tau(\x)-13,600 \le 0, \\ \vspace{3pt}
g_4(\x)=\sigma(\x)-30,000 \le 0, \\ \vspace{3pt}
g_5(\x)=0.10471 w^2 +0.04811  h d (14+L) -5.0 \le 0, \\ \vspace{3pt}
g_6(\x)=0.125 - w \le 0, \\ \vspace{3pt}
g_7(\x)=6000 - P(\x) \le 0,
\end{array}
\ee
where
\be \begin{array}{ll}
 \sigma(\x)=\frac{504,000}{h d^2},  & Q=6000 (14+\frac{L}{2}), \\ \\
 D=\frac{1}{2} \sqrt{L^2 + (w+d)^2}, & J=\sqrt{2} \; w L [ \frac{L^2}{6} + \frac{(w+d)^2}{2}], \\ \\
 \delta=\frac{65,856}{30,000 h d^3}, &  \beta=\frac{QD}{J}, \\ \\
 \alpha=\frac{6000}{\sqrt{2} w L}, & \tau(\x)=\sqrt{\alpha^2 + \frac{\alpha \beta L}{D}+\beta^2}, \\ \\
 P=0.61423 \times 10^6 \; \frac{d h^3}{6} (1-\frac{d \sqrt{30/48}}{28}). &
\end{array} \ee
The simple limits or bounds are $0.1 \le L, d \le 10$ and
$0.1 \le w, h \le 2.0$.

Using our Cuckoo Search, we have the following optimal solution
\[ \x_*=(w ,L, d, h) \] \be =
(0.205729639786079,   3.470488665627977,   9.036623910357633,   0.205729639786079), \ee
with \be f(\x*)_{\min} = 1.724852308597361.      \ee
This solution is exactly the same as the solution obtained by Cagnina et al (2008)
\be f_*=1.724852 \quad \textrm{ at } \;\; (0.205730, 3.470489, 9.036624, 0.205729). \ee
We have seen that, for both test problems, CS has found the optimal solutions which are
either better than or the same as the solutions found so far in literature.

\section{Discussions and Conclusions}

From the comparison study of the performance of CS with GAs and PSO, we know that
our new Cuckoo Search in combination with L\'evy flights is
very efficient and proves to be superior for almost all the test problems.
This is partly due to the fact that there are fewer parameters to be fine-tuned in
CS than in PSO and genetic algorithms. In fact, apart from the population size
$n$, there is essentially one parameter $p_a$. If we look at the CS algorithm carefully,
there are essentially three components: selection of the best, exploitation by local
random walk, and exploration by randomization via L\'evy flights globally.

The selection of the best by keeping the best nests or solutions
is equivalent to some form of elitism commonly used in genetic algorithms, which ensures the best
solution is passed onto the next iteration and there is no risk that the best
solutions are cast out of the population. The exploitation around the
best solutions is performed by using a local random walk
\be \x^{t+1}=\x^t + \a \ff{\varepsilon}_t. \ee
If $\ff{\varepsilon}_t$ obeys a Gaussian distribution, this becomes a standard random walk indeed.
This is equivalent to the crucial step in pitch adjustment
in Harmony Search (Geem et al 2001, Yang 2009). If $\ff{\varepsilon}_t$ is drawn from a L\'evy
distribution, the step of move is larger, and could be potentially more efficient. However,
if the step is too large, there is risk that the move is too far away. Fortunately, the
elitism by keeping the best solutions makes sure that the exploitation moves are within
the neighbourhood of the best solutions locally.

On the other hand, in order to sample the search space effectively so that new
solutions to be generated are diverse enough, the exploration step is carried out in terms
of L\'evy flights. In contrast, most metaheuristic algorithms use either uniform distributions
or Gaussian to generate new explorative moves (Geem et al 2001, Blum and Rilo 2003).
If the search space is large, L\'evy flights are usually more efficient.
A good combination of the above three components can thus lead to an efficient algorithm such as Cuckoo Search.

Furthermore, our simulations also
indicate that the convergence rate is insensitive to the algorithm-dependent
parameters such as $p_a$.
This also means that we do not have to fine tune these parameters for
a specific problem. Subsequently, CS is more generic and robust for
many optimisation problems, comparing with other metaheuristic algorithms.

This potentially powerful optimisation strategy can easily be extended
to study multiobjecitve optimization applications with various constraints,
including NP-hard problems. Further studies can focus  on the
sensitivity and parameter studies and their possible relationships
with the convergence rate of the algorithm. In addition, hybridization with other popular
algorithms such as PSO will also be potentially fruitful.
More importantly, as for most metaheuristic algorithms, mathematical analysis
of the algorithm structures is highly needed. At the moment, no
such framework exists for analyzing metaheuristics in general. Any progress in
this area will potentially provide new insight into the understanding of how and why
metaheuristic algorithms work.

\section*{Appendix: Demo Implementation}

{\small 
\begin{verbatim}
% -----------------------------------------------------------------
% Cuckoo Search (CS) algorithm by Xin-She Yang and Suash Deb      %
% Programmed by Xin-She Yang at Cambridge University              %
% Programming dates: Nov 2008 to June 2009                        %
% Last revised: Dec  2009   (simplified version for demo only)    %
% -----------------------------------------------------------------
% Papers -- Citation Details:
% 1) X.-S. Yang, S. Deb, Cuckoo search via Levy flights,
% in: Proc. of World Congress on Nature & Biologically Inspired
% Computing (NaBIC 2009), December 2009, India,
% IEEE Publications, USA,  pp. 210-214 (2009).
% http://arxiv.org/PS_cache/arxiv/pdf/1003/1003.1594v1.pdf 
% 2) X.-S. Yang, S. Deb, Engineering optimization by cuckoo search,
% Int. J. Mathematical Modelling and Numerical Optimisation, 
% Vol. 1, No. 4, 330-343 (2010). 
% http://arxiv.org/PS_cache/arxiv/pdf/1005/1005.2908v2.pdf
% ----------------------------------------------------------------%
% This demo program only implements a standard version of         %
% Cuckoo Search (CS), as the Levy flights and generation of       %
% new solutions may use slightly different methods.               %
% The pseudo code was given sequentially (select a cuckoo etc),   %
% but the implementation here uses Matlab's vector capability,    %
% which results in neater/better codes and shorter running time.  % 
% This implementation is different and more efficient than the    %
% the demo code provided in the book by 
%    "Yang X. S., Nature-Inspired Metaheuristic Algoirthms,       % 
%     2nd Edition, Luniver Press, (2010).                 "       %
% --------------------------------------------------------------- %

% =============================================================== %
% Notes:                                                          %
% Different implementations may lead to slightly different        %
% behavour and/or results, but there is nothing wrong with it,    %
% as this is the nature of random walks and all metaheuristics.   %
% -----------------------------------------------------------------

function [bestnest,fmin]=cuckoo_search(n)
if nargin<1,
% Number of nests (or different solutions)
n=25;
end

% Discovery rate of alien eggs/solutions
pa=0.25;

%% Change this if you want to get better results
% Tolerance
Tol=1.0e-5;
%% Simple bounds of the search domain
% Lower bounds
nd=15; 
Lb=-5*ones(1,nd); 
% Upper bounds
Ub=5*ones(1,nd);

% Random initial solutions
for i=1:n,
nest(i,:)=Lb+(Ub-Lb).*rand(size(Lb));
end

% Get the current best
fitness=10^10*ones(n,1);
[fmin,bestnest,nest,fitness]=get_best_nest(nest,nest,fitness);

N_iter=0;
%% Starting iterations
while (fmin>Tol),

    % Generate new solutions (but keep the current best)
     new_nest=get_cuckoos(nest,bestnest,Lb,Ub);   
     [fnew,best,nest,fitness]=get_best_nest(nest,new_nest,fitness);
    % Update the counter
      N_iter=N_iter+n; 
    % Discovery and randomization
      new_nest=empty_nests(nest,Lb,Ub,pa) ;
    
    % Evaluate this set of solutions
      [fnew,best,nest,fitness]=get_best_nest(nest,new_nest,fitness);
    % Update the counter again
      N_iter=N_iter+n;
    % Find the best objective so far  
    if fnew<fmin,
        fmin=fnew;
        bestnest=best;
    end
end %% End of iterations

%% Post-optimization processing
%% Display all the nests
disp(strcat('Total number of iterations=',num2str(N_iter)));
fmin
bestnest

%% --------------- All subfunctions are list below ------------------
%% Get cuckoos by ramdom walk
function nest=get_cuckoos(nest,best,Lb,Ub)
% Levy flights
n=size(nest,1);
% Levy exponent and coefficient
% For details, see equation (2.21), Page 16 (chapter 2) of the book
% X. S. Yang, Nature-Inspired Metaheuristic Algorithms, 2nd Edition, Luniver Press, (2010).
beta=3/2;
sigma=(gamma(1+beta)*sin(pi*beta/2)/(gamma((1+beta)/2)*beta*2^((beta-1)/2)))^(1/beta);

for j=1:n,
    s=nest(j,:);
    % This is a simple way of implementing Levy flights
    % For standard random walks, use step=1;
    %% Levy flights by Mantegna's algorithm
    u=randn(size(s))*sigma;
    v=randn(size(s));
    step=u./abs(v).^(1/beta);
  
    % In the next equation, the difference factor (s-best) means that 
    % when the solution is the best solution, it remains unchanged.     
    stepsize=0.01*step.*(s-best);
    % Here the factor 0.01 comes from the fact that L/100 should the typical
    % step size of walks/flights where L is the typical lenghtscale; 
    % otherwise, Levy flights may become too aggresive/efficient, 
    % which makes new solutions (even) jump out side of the design domain 
    % (and thus wasting evaluations).
    % Now the actual random walks or flights
    s=s+stepsize.*randn(size(s));
   % Apply simple bounds/limits
   nest(j,:)=simplebounds(s,Lb,Ub);
end

%% Find the current best nest
function [fmin,best,nest,fitness]=get_best_nest(nest,newnest,fitness)
% Evaluating all new solutions
for j=1:size(nest,1),
    fnew=fobj(newnest(j,:));
    if fnew<=fitness(j),
       fitness(j)=fnew;
       nest(j,:)=newnest(j,:);
    end
end
% Find the current best
[fmin,K]=min(fitness) ;
best=nest(K,:);

%% Replace some nests by constructing new solutions/nests
function new_nest=empty_nests(nest,Lb,Ub,pa)
% A fraction of worse nests are discovered with a probability pa
n=size(nest,1);
% Discovered or not -- a status vector
K=rand(size(nest))>pa;

% In the real world, if a cuckoo's egg is very similar to a host's eggs, then 
% this cuckoo's egg is less likely to be discovered, thus the fitness should 
% be related to the difference in solutions.  Therefore, it is a good idea 
% to do a random walk in a biased way with some random step sizes.  
%% New solution by biased/selective random walks
stepsize=rand*(nest(randperm(n),:)-nest(randperm(n),:));
new_nest=nest+stepsize.*K;

% Application of simple constraints
function s=simplebounds(s,Lb,Ub)
  % Apply the lower bound
  ns_tmp=s;
  I=ns_tmp<Lb;
  ns_tmp(I)=Lb(I);
  
  % Apply the upper bounds 
  J=ns_tmp>Ub;
  ns_tmp(J)=Ub(J);
  % Update this new move 
  s=ns_tmp;

%% You can replace the following by your own functions
% A d-dimensional objective function
function z=fobj(u)
%% d-dimensional sphere function sum_j=1^d (u_j-1)^2. 
%  with a minimum at (1,1, ...., 1); 
z=sum((u-1).^2);

\end{verbatim}

}


\begin{thebibliography}{100}

\bibitem{Arora}
Arora, J., 1989. {\it Introduction to Optimum Design}, McGraw-Hill.

\bibitem{Bele}
Belegundu, A., 1982. `A study of mathematical programming methods
for structural optimization', PhD thesis, Department of Civil Environmental
Engineering, University of Iowa, USA.

\bibitem{Bath}
Barthelemy, P., Bertolotti, J., Wiersma, D. S., 2008.
`A L\'evy flight for light', {\it Nature}, {\bf 453}, 495-498.

\bibitem{Blum}
Blum, C. and Roli, A., 2003. `Metaheuristics in combinatorial optimization: Overview and
conceptural comparision', {\it ACM Comput. Surv.}, {\bf 35}, 268-308.

\bibitem{Brown}
Brown, C., Liebovitch, L. S., Glendon, R., 2007. `L\'evy flights in Dobe Ju/'hoansi
foraging patterns', {\it Human Ecol.}, {\bf 35}, 129-138.

\bibitem{Cag}
Cagnina, L. C., Esquivel, S. C., and Coello, C. A., 2008.
`Solving engineering optimization problems with the simple
constrained particle swarm optimizer', {\it Informatica},
{\bf 32}, 319-326.


\bibitem{Chatt}
Chattopadhyay, R., 1971. `A study of test functions for optimization
algorithms', {\it J. Opt. Theory Appl.}, {\bf 8}, 231-236.


\bibitem{Deb}
Deb, K., 1995. {\it Optimisation for Engineering Design}, Prentice-Hall, New Delhi.


\bibitem{Floudas}
Floudas, C. A., Pardalos, P. M., Adjiman, C. S., Esposito, W. R.,
Gumus, Z. H., Harding, S. T.,  Klepeis, J. L., Meyer, C. A., Scheiger, C. A., 1999.
{\it Handbook of Test Problems in Local and Global Optimization}, Springer, 1999.


\bibitem{Gazi}
Gazi, K., and Passino, K. M., 2004. Stability analysis of social foraging swarms,
{\it IEEE Trans. Sys. Man. Cyber. Part B - Cybernetics}, {\bf 34}, 539-557.

\bibitem{Geem}
Geem, Z. W.,   Kim, J. H.,  Loganathan, G. V., 2001.
`A new heuristic optimization algorithm: Harmony search', {\it Simulation}, {\bf 76}, 60-68.



\bibitem{Gold}
Goldberg, D. E., 1989. {\it Genetic Algorithms in Search, Optimisation and
Machine Learning}, Reading, Mass., Addison Wesley.

\bibitem{Hedar}
Hedar, A., 2005, `Test function web pages',
http://www-optima.amp.i.kyoto-u.ac.jp /member/student/hedar/Hedar$\_$files/TestGO$\_$files/Page364.htm


\bibitem{Ken}
Kennedy, J. and Eberhart, R. C., 1995. `Particle swarm optimization'. {\it
Proc. of IEEE International Conference on Neural Networks},
Piscataway, NJ. pp. 1942-1948.

\bibitem{Ken2}
Kennedy, J., Eberhart, R. C., Shi, Y., 2001. {\it Swarm intelligence}, Academic Press.

\bibitem{Molga}
Molga, M., Smutnicki, C., 2005. ``Test functions for optimization needs'', \\
http://www.zsd.ict.pwr.wroc.pl/files/docs/functions.pdf


\bibitem{Passino}
Passino, K. M., 2001. {\it Biomimicry of Bacterial Foraging for Distributed
Optimization}, University Press, Princeton, New Jersey.

\bibitem{Payne}
Payne, R. B., Sorenson, M. D., and Klitz, K.,2005. {\it The Cuckoos},
Oxford University Press.

\bibitem{Pav}
Pavlyukevich, I., 2007. `L\'evy flights, non-local search and simulated annealing',
{\it J. Computational Physics}, {\bf 226}, 1830-1844.


\bibitem{Rag}
Ragsdell, K. and Phillips, D.,1976. `Optimal design of a class of welded
structures using geometric programming', {\it J. Eng. Ind.},
{\bf 98}, 1021-1025.

\bibitem{Rey}
Reynolds, A. M. and Frye, M. A., 2007. `Free-flight odor tracking in
Drosophila is consistent with an optimal intermittent scale-free search',
{\it PLoS One}, {\bf 2}, e354.

\bibitem{Schoen}
Schoen, F., 1993. `A wide class of test functions for
global optimization', {\it J. Global Optimization},
{\bf 3}, 133-137.


\bibitem{Shang}
Shang, Y. W., Qiu Y. H., 2006. `A note on the extended rosenrbock function',
{\it Evolutionary Computation}, {\bf 14}, 119-126.


\bibitem{Shilane}
Shilane D., Martikainen J., Dudoit S., Ovaska S. J., 2008.
 `A general framework for statistical performance comparison
of evolutionary computation algorithms', {\it Information Sciences},
{\bf 178}, 2870-2879.


\bibitem{Shles}
Shlesinger, M. F.,2006. `Search research', {\it Nature},
{\bf 443}, 281-282.


\bibitem{Yang}
Yang, X. S., 2008. {\it Nature-Inspired Metaheuristic Algorithms},
Luniver Press, (2008).

\bibitem{Yang2}
Yang, X. S., 2005. `Biology-derived algorithms in engineering optimizaton' (chapter 32),
in {\it Handbook of Bioinspired Algorithms and Applications}
(eds Olarius \& Zomaya),  Chapman \& Hall / CRC.

\bibitem{YangDeb}
Yang, X. S. and Deb, S., 2009. `Cuckoo search via L\'evy flights',
{\it Proceeings of World Congress on Nature \& Biologically Inspired
Computing} (NaBIC 2009, India), IEEE Publications, USA, pp. 210-214.

\bibitem{Yang4}
Yang, X. S., 2009. `Harmony search as a metaheuristic algorithm',
in: {\it Music-Inspired Harmony Search: Theory and Applications} (Eds Z. W. Geem),
Springer, pp.1-14.


\bibitem{Yang5}
 Yang, X. S., 2010. {\it Engineering Optimisation: An Introduction with
 Metaheuristic Applications}, John Wiley and Sons.


\end{thebibliography}
\end{document}